\PassOptionsToPackage{usenames, dvipsnames}{color}

\documentclass[12pt,letterpaper,english]{amsart}

\usepackage{amssymb,amsmath,amsthm,amscd,cite,mathrsfs,graphicx,multicol,booktabs,arydshln,stmaryrd}
\usepackage[cmtip,all,matrix,arrow,tips,curve]{xy}

\usepackage[margin=2.55cm]{geometry}
\usepackage{lipsum}

\usepackage{xypic}
\usepackage{tikz}
\usepackage{caption}
\usepackage{easybmat}
\usepackage{comment}
\newtheorem{thm}{Theorem}[section]

\newtheorem{lemma}[thm]{Lemma}

\theoremstyle{definition}
\newtheorem{remark}[thm]{Remark}

\newtheorem{question}[thm]{Question}



\newcommand\bH{{\mathbb H}}

\newcommand\bP{{\mathbb P}}

\newcommand\bV{{\mathbb V}}

\newcommand\cA{{\mathcal A}}
\newcommand\cB{{\mathcal B}}

\renewcommand\cD{{\mathcal D}}

\newcommand\cF{{\mathcal F}}

\renewcommand\cH{{\mathcal H}}

\renewcommand\cL{{\mathcal L}}

\newcommand\cM{{\mathcal M}}
\newcommand\cN{{\mathcal N}}
\newcommand\cO{{\mathcal O}}

\newcommand\cQ{{\mathcal Q}}

\newcommand\cS{{\mathcal S}}

\newcommand\cV{{\mathcal V}}

\newcommand\cX{{\mathcal X}}
\newcommand\cY{{\mathcal Y}}
\newcommand\cZ{{\mathcal Z}}

\newcommand\bfG{{\mathbf G}}
\newcommand\bfP{{\mathbf P}}

\newcommand{\nc}{\newcommand}
\nc{\rnc}{\renewcommand}
\rnc{\P}{\mathbb P}
\nc{\R}{\mathbb R}
\nc{\Q}{\mathbb Q}
\nc{\C}{\mathbb C}
\nc{\Z}{\mathbb Z}
\nc{\A}{\mathcal A}
\rnc{\H}{\mathbb H}
\rnc{\cO}{\mathcal O}
\rnc{\O}{\mathrm O}

\nc{\codim}{\mathrm{codim}}
\nc{\Sym}{\mathrm{Sym}}
\nc{\Herm}{\mathrm{Herm}}
\nc{\Spec}{\mathrm{Spec}\,}
\nc{\Proj}{\mathrm{Proj}\,}
\nc{\prim}{\mathrm{prim}}
\nc{\eps}{\epsilon}
\nc{\ov}{\overline}
\nc{\bs}{\backslash}

\nc{\M}{\mathcal M}
\nc{\D}{\mathbb D}
\nc{\B}{\mathbb B}
\rnc{\L}{\mathbb L}
\nc{\K}{\mathcal K}
\rnc{\S}{\mathcal S}
\nc{\Pic}{\mathrm{Pic}}
\nc{\Mat}{\mathrm{Mat}}
\nc{\Mod}{\mathrm{Mod}}
\nc{\QMod}{\mathrm{QMod}}
\nc{\SMod}{\mathrm{SMod}}
\nc{\Cl}{\mathrm{Cl}}
\nc{\End}{\mathrm{End}}
\nc{\NL}{\mathrm{NL}}
\nc{\Aut}{\mathrm{Aut}}
\nc{\rk}{\mathrm{rk}}
\nc{\E}{\mathrm{E}_8}
\nc{\U}{\mathrm{U}}
\nc{\SO}{\mathrm{SO}}
\nc{\GO}{\mathrm{GO}}
\nc{\Spin}{\mathrm{Spin}}
\nc{\Sp}{\mathrm{Sp}}
\nc{\SU}{\mathrm{SU}}
\nc{\SL}{\mathrm{SL}}
\nc{\GL}{\mathrm{GL}}
\nc{\CH}{\mathrm{CH}}
\nc{\Hilb}{\mathrm{Hilb}}
\nc{\vdim}{\mathrm{vdim}}
\nc{\pr}{\mathrm{pr}}
\nc{\mhm}{\mathrm{MHM}}
\nc{\Trace}{\mathrm{Tr}}
\nc{\cl}{\mathrm{cl}}
\nc{\vir}{\mathrm{vir}}

\theoremstyle{plain}

\newtheorem{theorem}{Theorem}[section]
\newtheorem{proposition}[theorem]{Proposition}

\newtheorem{conjecture}[theorem]{Conjecture}

\newtheorem{definition}[theorem]{Definition}

\theoremstyle{remark}

\newtheorem{example}[theorem]{Example}

\usepackage{cleveref}

\begin{document}

\pagenumbering{gobble}

\title{Modularity of special cycles on Shimura varieties:\\ a survey}
\author{Fran\c{c}ois Greer}
\address{Department of Mathematics, Michigan State University, 619 Red Cedar Rd, East Lansing, MI 48824}
\email{greerfra@msu.edu}
\author{Salim Tayou}
\address{Department of Mathematics, Dartmouth College, Kemeny Hall, Hanover, NH 03755, USA}
\email{salim.tayou@dartmouth.edu}

\begin{abstract}
We survey recent results on a conjecture of Kudla regarding the modularity of generating series of special cycle classes in toroidal compactifications of orthogonal and unitary Shimura varieties. Along the way, we formulate several conjectures on related phenomena for special cycles in other types of Shimura varieties, as well as on more general quotients of period domains.

\end{abstract}

\maketitle
\section{Introduction}

In their groundbreaking paper, Hirzebruch and Zagier \cite{hirzag} proved that the generating series of corrected special cycles on a compactified Hilbert modular surface is a modular form of weight $2$, with values in the cohomology. Their work prompted a series of developments culminating in the theorem of Kudla and Millson \cite{kudla-millson}, which gave a general modularity statement for special cycles on orthogonal and unitary Shimura varieties. Our goal in this note is twofold: to propose an analogue of the Kudla-Millson theorem for other classes of Shimura varieties, and to formulate a conjectural framework for how the special cycles can be completed and corrected in toroidal compactifications of these Shimura varieties.

\medskip 

Let $(\bfG,\cD)$ be a Shimura datum; that is, $\bfG$ is a semisimple algebraic group over $\Q$, and $\cD$ is a Hermitian symmetric domain equipped with a transitive action of $G= \bfG(\R)$:
\[G\rightarrow \mathrm{Aut}(\cD), \]with compact stabilizer $K$, satisfying the axioms of Deligne \cite[Section 1.5]{deligneshimura}. Of particular interest are the irreducible Hermitian symmetric domains $\cD= G/K$ which have been completely classified: they consist of four infinite families (Types I--IV), together with two exceptional cases ($E_6$ and $E_7$). We will focus only on the infinite families here.

\begin{table}[!h]
\begin{center}
\begin{tabular}{|c|c|c|c|c|c|} 
 \hline
Type & $G$ & $\bfG_\C$ & $K$ &  $\dim_\C(X)$ & moduli problem  \\ \hline 
I & $\SU(a,b)$ & $\SL_{a+b}$  & $S(\U(a)\times \U(b))$ & $ab$ & pol. a.v. with CM \\ \hline
II & $\SO^*(2n)$ & $\SO_{2n}$ & $\U(n)$ & $n(n-1)/2$ & pol. a.v. with QM  \\ \hline
III & $\Sp(2n,\R)$ & $\Sp_{2n}$ & $\U(n)$ & $n(n+1)/2$ & pol. a.v. \\ \hline
IV & $\SO(n,2)$ & $\SO_{n+2}$ & $S(\O(n)\times \O(2))$ & $n$ & Kuga-Satake a.v.\\
\hline
\end{tabular}\bigskip
\caption{The four infinite families, with $n>0$ and $a\geq b>0$.}
\label{table_sv}
\end{center}
\end{table}


The choice of a faithful representation $\bfG\rightarrow \mathrm{GL}(V_\Q)$, where $V$ is a free $\Z$-module, determines an arithmetic subgroup $\Gamma$ of $\bfG(\Q)$ consisting of the elements that preserve the lattice $V$. The double quotient
$$X = X_\Gamma = \Gamma \bs G/K$$
is a complex quasi-projective variety by the theorem of Baily-Borel \cite{bailyborel}.
\medskip

There are natural families of abelian varieties $\mathcal B$ with polarization, endomorphism algebra, and level structure (PEL) over the Shimura varieties of Types I--III. In each of these cases, $X$ supports a PVHS of weight 1, where the fiber at each point $x\in X$ is isomorphic to
$V \simeq H^1(B_x,\Z).$
\medskip

The Type IV case is more subtle; one needs to take the double cover of $\SO(n,2)$, the spin group $\Spin(n,2)$ which admits a representation into the linear group of $\Cl(V)$, the Clifford algebra of $V$. The Spin Shimura variety admits a family $\mathcal B$ of polarized abelian varieties of dimension $2^{n+1}$ given by the Kuga-Satake construction. We have an isomorphism:

$$\Cl(V)\simeq H^1(B_x,\Z),$$
and an embedding $V\hookrightarrow \End(H^1(B_x,\Z))$.

For each of the types I--IV, we will recall the construction of certain natural classes of special cycles on $X$, which are always effective, and how to assemble them into generating series. Some of these constructions go back to the original work of Kudla, Millson, Zagier, and Hirzebruch, while others have only recently been considered by various groups of authors. The first goal of this survey is to place these constructions into a unified framework.
\medskip 

The simplest class of special cycles are indexed by integers $m\in \mathbb N$, and we denote them by $\cZ(m)$. The classes $[\cZ(m)]$ in the relevant cohomology or Chow group of $X$ constitute\footnote{This is a theorem in Types (I, IV) and a conjecture in Types (II, III)} the Fourier coefficients of a classical modular form of weight $k$ (recorded below) and level determined by the arithmetic group $\Gamma$. 
$$\Phi(q) = \sum_{m\geq 0} [\cZ(m)]q^m \in \CH^{c}(X) \otimes_\Q \Mod_k(\SL_2(\Z))$$
The constant term $[\cZ(0)]$ is defined separately as the top Chern class of a tautological bundle on $X$, so it only makes sense as a cycle class, and it is not effective in general.

\begin{table}[!h]
\begin{center}
\begin{tabular}{|c|c|c|c|c|c|} 
 \hline
Type & $G$ & $\dim_\C(X)$ & $c=\codim(Z(m))$ &  weight $k$ & moduli description  \\ \hline 
I & $\SU(a,b)$ & $ab$  & $a$ & $a+b$ & Noether--Lefschetz \\ \hline
II & $\SO^*(2n)$ & $n(n-1)/2$ & $n-1$ & $2n$ & isogeny split \\ \hline
III & $\Sp(2n,\R)$ & $n(n+1)/2$ & $n-1$ & $2n$ & isogeny split \\ \hline
IV & $\SO(n,2)$ & $n$ & $1$ & $(n+2)/2$ & Noether--Lefschetz\\
\hline
\end{tabular}\bigskip
\caption{Special cycles on Shimura varieties}
\label{table_cycles}
\end{center}
\end{table}

As we will explain, the definitions of these cycles are qualitatively different between Types (I, IV) and Types (II, III). For the former two types, each $\cZ(m)$ is a union of sub-Shimura varieties of that same type, whereas for the latter two types, each $\cZ(m)$ is a Hecke translate of a particular sub-Shimura variety $\cZ(1)$ that comes from a product of two irreducible Hermitian symmetric domains of that same type.

For each fixed $g\geq 1$, one can generalize the definitions above to a class of special cycles $\cZ(M)$ of codimension $gc$, in Types (I, IV), indexed now by $M$ an integral positive definite symmetric (or Hermitian) matrix. Their classes $[\cZ(M)]$ constitute the Fourier coefficients of a Siegel modular form of the same weight $k$ from Table \ref{table_cycles}.
$$\Phi^g(q) = \sum_{M\geq 0} [\cZ(M)] q^M \in \CH^{gc}(X) \otimes_\Q \Mod_k(\Sp_{2g}(\Z))$$
The expressions $[\cZ(M)]$ when $M$ is only semi-positive are interpreted using an appropriate power of the tautological top Chern class $[\cZ(0)]$; see Section \ref{s:orthogonal}. In Section \ref{s:pel-case}, we will propose a conjectural Siegel modular series of higher codimension cycles in Type III.

\subsection{Compactifications}

For each Shimura variety $X$, the canonical  Satake-Baily-Borel compactification $X^*$ is a projective algebraic variety whose boundary is stratified by Shimura varieties of lower dimensions. Each boundary stratum $X'\subset X^*$ corresponds to a unique $\Gamma$-conjugacy class of maximal parabolic subgroups $\bfP\subset \bfG$ defined over $\Q$. The drawback of $X^*$ is that it is singular in general. A nice class of smooth compactifications was constructed in \cite{amrt}; these are called the \emph{toroidal compactifications} $X^\Sigma$, obtained by blowing up the boundary of $X^*$. They depend on a choice of admissible fan $\Sigma$ inside a certain cone $\mathscr C$ associated to the parabolic $\bfP$. Let
$$\eps: X^\Sigma \to X^*$$
denote this blow up map. Over each stratum $X'$ in the boundary, $\eps$ restricts to a proper surjective morphism which factors as
$$\eps^{-1}(X') \overset{\tau}\to \mathcal Y \overset{\pi}\to X',$$
where $\tau$ is a toric fibration, and $\pi$ is an abelian fibration. The second goal of this note is to explore the structure of the boundary of these compactifications in the cases of Table \ref{table_sv}, and to describe the limits of the special cycles from Table \ref{table_cycles} in the boundary of $X^\Sigma$. 
The Zariski closures $\overline{\cZ}(m)\subset X^{\Sigma}$ can be assembled into a generating series valued in the relevant cohomology or Chow group of $X^{\Sigma}$, but they are no longer modular forms; they are expected to be \emph{mixed mock modular forms} in general. Kudla conjectured that in types (I, IV), the modularity property can be restored after subtracting a series of effective cycle classes supported in the boundary of $X^\Sigma$.


\medskip

Special cycles on Shimura varieties are examples of Hodge loci and it is natural to explore generating series in the more general setting of polarized variations of (mixed) Hodge structure over an arbitrary base, which we explore in \Cref{s:vhs}. 

This survey is by no means exhaustive, and the interested reader should consult the original papers cited here as well as the excellent surveys \cite{kudla_special_2004}, \cite{Chao-li-ihes-notes}, \cite{chao-li-SW-formulas}. One can formulate everything we have written here on integral models $\mathcal X$ of Shimura varieties, where the corresponding modular forms are valued in arithmetic Chow groups of $\mathcal X$ and its compactifications. While this direction has been explored in detail in the literature for Types I and IV, we will say little about it in this article.

\subsection{Organization of the paper}
We survey classical constructions and results for orthogonal (IV) and unitary (I) Shimura varieties in \Cref{s:ortho-unitary-case}. In \Cref{s:pel-case}, we define special cycles in the more general setting of PEL Shimura varieties, which include Siegel modular varieties (III) and quaternionic Shimura varieties (II). In \Cref{s:compact}, we address the problem of extending the modularity results, first to Baily-Borel compactifications and then to toroidal compactifications in general. In \Cref{s:vhs}, we discuss the more general case of arithmetic quotients of non-classical period domains, which are not algebraic. Finally, in \Cref{s:coincidence} we speculate about possible relationships between the generating series for pairs of isogenous groups from different families in Table \ref{table_sv}.


\subsection{Acknowledgments} S.T. would like to thank the Institute for Advanced Study for its hospitality. S.T. is supported by NSF Award DMS-2503815 and by the Burke Research Initiation Award from Dartmouth College. F.G. is supported by NSF Award DMS-2302548. We would like to thank Carl Lian for helpful comments on an earlier draft.

\section{Orthogonal and unitary Shimura varieties}\label{s:ortho-unitary-case}
We recall here the construction of special cycles in orthogonal and unitary Shimura varieties following \cite{kudla-millson}, as well as Kudla's  modularity conjecture \cite{kudla_special_2004}.
\subsection{The orthogonal case}\label{s:orthogonal} Let $(L,Q)$ be a quadratic lattice of signature $(n,2)$ with $n\geq 1$. We denote by $(\cdot,\cdot)$ the symmetric bilinear form associated to $Q$ by the relation $Q(x)=\frac{(x\cdot x)}{2}$. To simplify the discussion, we assume that the lattice $(L,Q)$ is unimodular, which implies that the resulting modular forms will have level 1.
\medskip 

Let $\bfG$ be the orthogonal group of $(L_\Q,Q)$ and consider the Hermitian symmetric domain 

\[\cD=\{x\in \bP(L_{\C})\,\big|\, (x\cdot x)=0, (x\cdot\overline{x})<0\},\]
which parametrizes polarized K3-type Hodge structures on $L$. Then $(\bfG,\cD)$ is a Shimura datum of Type $IV$. The choice of an arithmetic subgroup $\Gamma\subset\bfG(\Q)$ preserving $L$ yields the \emph{orthogonal Shimura variety} $X=\Gamma\backslash \cD$, which is a quasi-projective algebraic variety. 
\medskip

Let $1\leq g\leq n$. For each $M\in\mathrm{Sym}_{g\times g}(\Z)_{\geq 0}$ positive semi-definite, consider the cycle:
\[\cZ(M)^{\textrm{naive}}=\Gamma\backslash\left(\bigcup_{\underset{\bigl((\lambda_i\cdot\lambda_j)\bigr)_{1\leq i,j\leq g}=M}{\lambda_1,\cdots,\lambda_g\in L}}\{x\in \cD\,|\, (x\cdot\lambda_i)=0,\,1\leq i\leq g\}\right)~.\]
This cycle has codimension $r(M)$, the rank of the matrix $M$. To obtain a cycle class in codimension $g$, we define: 
\[[\cZ(M)]:=[\cZ(M)^{\textrm{naive}}]\cap[\cL^\vee]^{g-r(M)}\in \mathrm{CH}^g(X)\]
where $\cL\rightarrow X$ denotes the tautological line bundle. We can form the generating series: 
\[\Theta^g_L(\tau)=\sum_{M\geq 0}[\cZ(M)]q^M\in \mathrm{CH}^g(X)\llbracket q \rrbracket,\]
where $q^M=e^{2\pi i\Trace(\tau M)}$ and $\tau\in \bH_g$ is an element of the Siegel upper-half space of genus $g$. The following theorem summarizes decades of active research. 

\begin{theorem}
    The generating series $\Theta^g_L$ is the $q$-expansion of a holomorphic modular form of genus $g$ and weight $\frac{n+2}{2}$ with respect to the Siegel modular group $\mathrm{Sp}_{2g}(\Z)$, that is:
    \[ \Theta^g_L\in \mathrm{CH}^g(X)\otimes_\Q {\mathrm{Mod}}_{\frac{n+2}{2}}(\mathrm{Sp}_{2g}(\Z))~.\]
\end{theorem}

The theorem above is a consequence of the work of many people: first in cohomology by Kudla-Millson \cite{kudla-millson}, then in $\mathrm{CH}^1$ by  Borcherds  \cite{borcherds-inventiones-automorphic}, in higher codimension by Wei Zhang \cite{zhang-thesis-2009} and Bruinier-Raum \cite{Bruinier-Westerholt-Raum}; see also \cite{howardmadapusi,howard-madapusi-kudla} for a different proof and an extension to integral models. 

\medskip

The Shimura variety $X$ is often non-compact; in fact it must be as soon as $n\geq 3$. It is natural to ask how the generating series $\Theta^g_L$ can be extended to a toroidal compactification of $X$, whence the following conjecture, due to Kudla \cite{kudla_special_2004}. 

\begin{conjecture}
    There exists a smooth toroidal compactification $X\hookrightarrow \overline{X}$, a finite collection $(Y_i)_{i\in I}$ of codimension $g$ cycles supported in the boundary $\overline X\smallsetminus X$, and corrected special cycles $\widetilde{\cZ}(M)=\overline{\cZ}(M)+\sum_ia_i(M)Y_i$ for $a_i(M)\in \Q$, such that the generating series: \[\sum_{M\geq 0}[\widetilde{\cZ}(M)]q^M\in\mathrm{CH}^g\left(\overline X\right)\llbracket q \rrbracket\]
is a holomorphic Siegel modular form of weight $\frac{n+2}{2}$ with respect to $\mathrm{Sp}_{2g}(\Z)$.
\end{conjecture}

This conjecture is still open in general, but the case of divisors was solved recently \cite{bruinierzemel, garcia-modularity,engel-greer-tayou}, as well as the case of zero-cycles \cite{bruinier-rosu-zemel}. We will give more details on these results in \Cref{s:compact}. 

\subsection{The unitary case}\label{s:unitary-sv}

Let $k$ be an imaginary quadratic field with ring of integers $\cO_k$, and let $(L,h)$ be a self-dual Hermitian $\cO_k$-module of signature $(a,b)$ with $a\geq b$. We fix an embedding $\tau:k\hookrightarrow \C$, which specifies a complex structure $J_0$ on $V_\R=V\otimes_\Q \R$.
\medskip 

The Hermitian symmetric domain associated to this data is the set:
\[\cD=\{P\subset L_\R\textrm{ of dimension } 2b\, \big|\, P \textrm{ is } J_0\textrm{-stable}\,\textrm{ and}\, h|_P \, \textrm{is negative definite}\}.\]

The group $\Gamma=\U(L)$ acts on $\cD$ and the quotient $X=\Gamma\backslash \cD$ is a unitary Shimura variety. 
\medskip 

Let $1\leq g\leq a$ and let $\mathrm{Herm}_g(\cO_k)$ denote the space of semi-positive Hermitian matrices with coefficients in $\cO_k$. For each $M\in \mathrm{Herm}_g(\cO_k)$, we have a special cycle \[\cZ(M)^{\textrm{naive}}=\Gamma\backslash\left(\bigcup_{\underset{\bigl((\lambda_i\cdot\lambda_j)\bigr)_{1\leq i,j\leq g}=M}{\lambda_1,\cdots,\lambda_g\in L}}\{x\in \cD\,|\, (x\cdot\lambda_i)=0,\,1\leq i\leq g\}\right)~.\]

This cycle has codimension $r(M) b$, where $r(M)$ is the rank of the matrix $M$. Let $\cV$ be the descent to $X$ of the tautological complex vector bundle on $\cD$ of rank $b$, and let $c_b(\cV)$ denote its top Chern class. To obtain a cycle class in codimension $gb$, we define 
\[[\cZ(M)]:=[\cZ(M)^{\textrm{naive}}]\cap [c_b(\cL^\vee)]^{g-r(M)}\in\mathrm{CH}^{gb}(X)~.\]

Let $\cH_g$ denote the Hermitian upper half-space and let $\U_{g,g}(\Z)$ be the integral Hermitian modular group. The special cycle classes defined above fit into the following generating series: 
\[\Theta^g_L(\tau)=\sum_{M\in \mathrm{Herm}_g(\cO_k)}[\cZ(M)] q^M\in\mathrm{CH}^{gb}(X)\llbracket q\rrbracket\]
 where $q^M=e^{2i\pi \Trace(\tau M)}$ and $\tau\in \cH_g$. The following result is due to Kudla and Millson \cite{kudla-millson}:
\begin{theorem}
    The generating series of cohomology classes $[\Theta^g_L]$ is the $q$-expansion of a holomorphic Hermitian modular form of weight $a+b$ with respect to $\U_{g,g}(\Z)$, that is:
    $$[\Theta_L^g] \in H^{2gb}(X,\Q) \otimes_\Q \Mod_{a+b}(\U_{g,g}(\Z)).$$
\end{theorem}
We note here that the modularity of the generating series $\Theta^g_L$ valued in the Chow group is still a major open problem; see Xia \cite{xia} for partial progress.

\medskip 

The Shimura variety $X$ must be non-compact as soon as $a+b\geq 3$, by Meyer's theorem, so one would like to understand the modularity behavior of the closures of special cycles in smooth compactifications. This motivates the following conjecture \cite[Conjecture 1.1]{bruinier-rosu-zemel}.
\begin{conjecture}
    There exists a smooth toroidal compactification $X\hookrightarrow \overline X$, a finite collection $(Y_i)_{i\in I}$ of codimension $gb$ cycles supported in the boundary $\overline X\smallsetminus X$, and corrected special cycles $\widetilde{\cZ}(M)=\overline{\cZ}(M)+\sum_ia_i(M)Y_i$ for $a_i(M)\in \Q$, such that the generating series: \[\sum_{M\geq 0}[\widetilde{\cZ}(M)]q^M\in\mathrm{CH}^{gb}(\overline X)\llbracket q\rrbracket\] is a holomorphic Hermitian modular form of weight $a+b$ with respect to $\mathrm{U}_{g,g}(\Z)$.
\end{conjecture}

For $g=b=1$, this conjecture was settled in \cite{bruinier-howard-kudla-yang}, and for $g=a$, $b=1$, it follows from the recent work 
\cite{bruinier-rosu-zemel}. In \Cref{s:compact}, we will discuss recent progress on the cohomological version of this conjecture for $b=1$, in joint work of the authors \cite{greer-tayou}.

\section{PEL Shimura varieties}\label{s:pel-case}

In this section, we formulate several conjectures on the modularity of certain special cycles in other Shimura varieties, some of which recover statements that have appeared in the literature. First, we deal with Siegel modular varieties and quaternionic Shimura varieties, and then we formulate analogous conjectures for general PEL Shimura varieties. A common theme in this section will be the presence of excess intersection classes, products of tautological Chern classes $\lambda_i = c_i(\mathbb E)$ which must be multiplied by an entire generating series of special cycles before formulating the precise modularity conjecture.

\subsection{Siegel modular varieties}\label{s:siegel-sv}
Let $(L,\psi)$ be the standard symplectic lattice of rank $2n$. Let $\bH_n$ denote the Siegel upper half-space. The quotient $\cA_n=\mathrm{Sp}_{2n}(\Z)\backslash \bH_n$ is a Siegel modular variety, interpreted as the moduli space of principally polarized abelian varieties. 
\medskip 

Using this modular interpretation, we define a class of special subvarieties in $\cA_n(\C)$ as follows. For each $m\geq 1$, let 
\[\cZ(m)=\{A\in \cA_n(\C)\,|\, \exists f:E\rightarrow A,\, f^\dagger\circ f= [m]_E, \, E\text{ elliptic curve}\}~.\]

Observe that the union $\bigcup_{m\geq 1}\cZ(m)$ corresponds exactly to the locus of abelian varieties admitting an elliptic curve as an isogeny factor. To form an algebraic cycle supported on each of these subvarieties $\cZ(m)$, we assign a finite multiplicity to each component equal to the {\it number} of maps $f:E \to A$ with norm $m$, for generic $A$ in that component of $\cZ(m)$. The cycles $\cZ(m)$ are in fact Hecke translates of $\cA_1\times \cA_{n-1}$, and hence they all have codimension $n-1$ inside $\cA_n$. Consider now the following generating series: 

\[\Theta_{\cA_n}(q)=\frac{(-1)^n}{24}\lambda_{n-1}+\sum_{m\geq 1}[\cZ(m)]q^m\in \mathrm{CH}^{n-1}(\cA_n)\llbracket q \rrbracket~,\]
The following conjecture appears in \cite[Conjecture 1]{greer-lian}.
\begin{conjecture}\label{conj_An}
    The generating series $\Theta_{\cA_n}(q)$ is a holomorphic modular form of weight $2n$ with respect to $\mathrm{SL}_2(\Z)$, valued in $\CH^{n-1}(\cA_n)$. 
\end{conjecture}

Partial progress on the cohomological version appears in \cite[Theorem 1.2]{greer-lian} with level structure, as well as in \cite[Corollary 6]{iribar-lopez}, where it is proved for $g\leq 5$. Using admissible covers of elliptic curves, \cite{lian} proved a similar statement for $g\leq 3$. 

A more systematic way to account for the multiplicities in the cycles $[\cZ(m)]$ is to introduce a moduli stack of triples:
$$\widetilde{\cZ}(m) = \{(E,A,f)\,|\, f:E \to A,\,f^\dagger\circ f= [m]_E\},$$
which has natural forgetful morphisms
$$\widetilde{\cZ}(m) \to \cA_1 \times \cA_n \overset{\pr}\longrightarrow \cA_n~.$$
Since $\cA_1$ is not proper, one cannot immediately pushforward the cycle classes $\left[\widetilde{\cZ}(m)\right]\in \CH^n(\cA_1\times \cA_n)$, but one does recover $[\cZ(m)]$ after first taking the closure in $\cA_1^*\times \cA_n$:
\[[\cZ(m)] = \pr_* \left[\cl\left(\widetilde{\cZ}(m)\right)\right] .\] 

To generalize Conjecture \ref{conj_An} to higher codimension, fix an integer $g$ with $1< g\leq n$. For any positive definite $M\in\mathrm{Sym}_{g\times g}(\Z)$, we define the following special subvariety of $\cA_n(\C)$:
\[\cZ(M)=\{A\in \cA_n(\C)\,|\, \exists f:E^g\rightarrow A,\, f^\dagger\circ f= [M]_{E^g}, \, E\text{ elliptic curve}\}~.\]
To form an algebraic cycle supported on each of these subvarieties $\cZ(M)$, we again weight each component by the generic number of such maps $f:E^g \to A$:
\[[\cZ(M)] = \pr_* \left[\cl\left(\widetilde{\cZ}(M)\right)\right] .\]

The resulting cycles have codimension $gn-1-{g \choose 2}$ inside $\cA_n$, since they are Hecke translates of the image of the finite map $\cA_1 \times \cA_{n-g} \to \cA_n$ sending $(E,B)\mapsto E^g\times B$. We form the generating series of special cycle classes:
$$\Theta_{\cA_n}^g(\tau) = \sum_{M\geq 0} [\cZ(M)] q^M \in \CH^{gn-1-{g\choose 2}}(\cA_n) \llbracket q \rrbracket~,$$
where $q^M=e^{2\pi i \Trace(\tau M)}$ and $\tau\in \mathbb H_g$, as in \Cref{s:ortho-unitary-case}. 

\begin{question}\label{Q:siegel}
    What can be said about Siegel modularity of the series $\Theta_{\cA_n}^g$?
\end{question}

We believe that a related modularity statement is true for {\it virtual} special cycles, defined by multiplying the normalized cycles $\left[\cl\left(\widetilde{\cZ}(M)\right)\right]$ with an excess intersection term $\lambda_1^{{g\choose 2}}$ pulled back from the factor $\cA_1^*$, before applying $\pr_*$. Consider the virtual series:
$$\Theta^{g,\vir}_{\cA_n}(\tau) = \sum_{M\geq 0} 
\pr_*\left[\lambda_1^{g\choose 2}\cdot \cl\left(\widetilde{\cZ}(M)\right)\right] q^M \in \CH^{gn-1}(\cA_n)\llbracket q \rrbracket~.$$
Observe that since $\dim(\cA_1^*)=1$, these  Fourier coefficients are nonzero only if $g\leq 2$.

\begin{conjecture}
The series $\Theta^{2,\vir}_{\A_n}$ is a holomorphic Siegel modular form of weight $2n$ with respect to $\Sp_{4}(\Z)$ valued in $\CH^{2n-1}(\cA_n)$, which vanishes at all cusps.
\end{conjecture}

Generalizing in a different direction, fix now $1< r \leq n/2$. For each $m\geq 1$, define the following special subvariety of $\cA_n(\C)$: 
\[\cZ_{r}(m)=\{A\in \cA_n\,|\,\exists f:B\rightarrow A, f^\dagger\circ f=[m]_B,\, B \textrm{ ppav of dimension } r\}~.\]

The union of the loci $\cZ_r(m)$ over $m\geq 1$ corresponds to the locus of principally polarized abelian varieties of dimension $g$ admitting an $r$-dimensional principally polarized abelian variety as an isogeny factor. They are Hecke translates of the cycles $\cA_r\times\cA_{n-r}$, and hence they all have codimension equal to $\binom{n+1}{2}-\binom{r
+1}{2}-\binom{n-r+1}{2}=r(n-r)$. We will form an algebraic cycle supported on each $\cZ_r(m)$ by assigning multiplicities in the usual way:
\[\widetilde{\cZ}_r(m) = \{ (B,A,f)\,|\, f:B\rightarrow A, f^\dagger\circ f=[m]_B  \} \to \cA_r\times \cA_n \overset{\pr}\longrightarrow \cA_n~.\]
Next, we form the generating series of cycle classes $[\cZ_r(m)]\in \CH^{r(n-r)}(\cA_n)$, defined now using the closure of $\widetilde{\cZ}_r(m)$ in $\cA_h'\times \cA_n$, for any smooth toroidal compactification $\cA_h\subset \cA_h'$:
\[[\cZ_r(m)] = \pr_* \left[\cl\left(\widetilde{\cZ}_r(m)\right)\right] .\] 
\[\Theta_{\mathrm{Sp}_{2n}}^r(q)=\sum_{m\geq 0}[\cZ_r(m)] q^m\in\mathrm{CH}^{r(n-r)}(\cA_n)\llbracket q\rrbracket~.\]
\begin{question}\label{Q:siegel2}
    What can be said about the modularity of the series $\Theta_{\mathrm{Sp}_{2n}}^r(q)$?
\end{question}
To provide some motivation for the questions above: let $(L_r,\psi)$ be the standard symplectic lattice of rank $2r$. The tensor product of lattices $(L,Q)=(L,\psi)\otimes (L_r,\psi_r)$ is endowed with a quadratic form of signature $(2rn,2rn)$. Let $Y$ be the arithmetic quotient of the period domain that parametrizes polarized Hodge structures on $(L,Q)$. We get an embedding: \[\cA_r\times\cA_n\hookrightarrow Y~.\]
A natural approach to Questions \ref{Q:siegel} and \ref{Q:siegel2} would be to pull back the orthogonal type Kudla--Millson generating series from $Y$. Several difficulties appear due to excess intersections when $r>1$, which motivate the following definition of the virtual generating series:
$$\Theta^{r,\vir}_{\Sp_{2n}}(q) =  \sum_{m\geq 0} 
\pr_*\left[\lambda_1\lambda_2\dots\lambda_{r-1}\cdot \cl\left(\widetilde{\cZ}_r(m)\right)\right] q^m \in \CH^{r(n-r)+{r\choose 2}}(\A_n)\llbracket q \rrbracket~. $$

\begin{conjecture}
    The series $\Theta^{r,\vir}_{\Sp_{2n}}(q)$ is a holomorphic modular form of weight $2gh$ with respect to $\SL_2(\Z)$ valued in $\CH^{r(n-r)+{r\choose 2}}(\cA_n)$.
\end{conjecture}

\subsection{Quaternionic Shimura varieties}\label{s:quaternion-SV}

The final class of Shimura varieties that we investigate in this section are those of Type II. 
\medskip 

Let $B$ be a quaternion algebra over $\Q$. We denote by $\alpha\rightarrow \alpha^*$ the main involution on $B$ and let $V$ be a finite-dimensional vector space over $\Q$. 
\begin{definition}
We say that $V$ is a quaternionic skew-Hermitian space if it admits a right $B$-action and a pairing \[(\cdot,\cdot):V\times V\rightarrow B\] such that 
\[\forall v,w\in V, \alpha,\beta\in B:\, (v\alpha,w\beta)=\alpha^*(v,w)\beta,\quad (v,w)=-(w,b)^*~.\]
\end{definition}

Let $V$ be a quaternionic skew-Hermitian space. The unitary similitude group, denoted $\bfG=\mathrm{GU}_B(V)$, is a algebraic group over $\Q$ which satisfies for each $\Q$-algebra $R$ 
\[\mathrm{GU}(R)=\{g\in \GL(V\otimes_\Q R)|\, (g(v),g(w))=\nu(g)(v,w),\, \forall v,w, \textrm{and some}\, \nu(g)\in R^{\times}\}~.\]

We assume henceforth that $B$ is ramified over $\R$. It follows that $\bfG_\R\simeq \mathrm{GO}^{*}(2n)$, where $n=\dim_B(V)$ is the rank of $V$ as a $B$-module. The Hermitian symmetric space associated to $\mathrm{GO}^{*}(2n)$ is given by: 
\[\mathcal H_{SO^*_{2n}}=\{Z\in M_n(\C)\,|\, ^tZZ+1=0, \mathrm{Im}(Z)>0\}~.\]

Let $\mathcal{O}$ be an order in $B$ which is stable under the main involution of $B$ and let $L\subset V$ be a finite rank $\Z$-module which is $\mathcal{O}$-stable. Let $\Gamma$ be the stabilizer of $L$ in $\bfG(\Q)$. Then the quotient \[\cQ_{n}:=\Gamma\backslash \cH_{\SO^{*}_{2n}}\] is a \emph{quaternionic Shimura variety.} 

\medskip

In order to define special cycles in $\Gamma\backslash \cH_{\SO^{*}_{2n}}$ with good moduli interpretation, we will follow an idea similar to that of the Siegel modular variety discussed in \Cref{s:siegel-sv}. Consider a $B$-orthogonal splitting of $V=V_1\oplus V_2$ such that the restriction of the skew Hermitian form is non-degenerate on each factor. This decomposition yields an embedding of algebraic groups \[\GO^*(2r)\times \GO^*(2n-2r)\hookrightarrow \GO^*(2n)~,\] as well as an embedding of Hermitian symmetric spaces: 
\[\cH_{\SO^{*}_{2r}}\times\cH_{\SO^{*}_{2n-2r}}\hookrightarrow\cH_{\SO^{*}_{2n}}~.\]
\medskip

Let $L_1\subseteq L\cap V_1$ and $L_2\subseteq L\cap V_2$ be finite index sub $\cO$-modules. Then $L_1\oplus L_2$ is a finite index sub $\cO$-module of $L$ and we denote its index by $m$. 

Let $\Gamma_1$ be the stabilizer of $L_1$ and $\Gamma_2$ the stabilizer of $L_2$. Then we get a finite morphism of Shimura varieties: 
\[\iota_{L_1,L_2}:\Gamma_1\backslash \cH_{\SO^{*}_{2r}}\times \Gamma_2\backslash \cH_{\SO^{*}_{2n-2r}}\rightarrow \cQ_{n}~.\]

For $m\geq 1$, let $\cZ(m)$ denote the union of the images \[\iota_{L_1,L_2}\left(\Gamma\backslash \cH_{\SO^{*}_{2r}}\times \Gamma\backslash \cH_{\SO^{*}_{2n-2r}}\right)\] over all possible decompositions $L_1\oplus L_2$ of index $m$ in $L$. By well-known results of Borel and Harish-Chandra, $\cZ(m)$ is an algebraic subvariety of $\cQ_{n}$ and it has codimension $\binom{n}{2}-\binom{r}{2}-\binom{n-r}{2}=r(n-r)$. We can now form the generating series: 
\[\Theta^r_{\cQ_n}(q)=\sum_{m\geq 1}[\cZ(m)]q^m\in \CH^{r(n-r)}(\cQ_{n})\llbracket q\rrbracket~.\]

\begin{question}
    What can be said about the modularity of the generating series $\Theta^k_{\cQ_n}$?
\end{question}

\begin{remark}
The reader will notice the similarities with Question \ref{Q:siegel2} and how one can get similar heuristics for the modularity statement above. 
    
\end{remark}

In the case where the Shimura variety $\cQ_{n}$ is non-compact, one can pose a conjecture similar to Kudla's about toroidal corrections.

\begin{question}
    Does there exist a toroidal compactification of $\cQ_{n}$ and correction of the cycles $\cZ(m)$ by boundary cycles such that resulting generating series is again a modular form of the same weight and level? 
\end{question}

\subsection{PEL Shimura varieties} 

The Shimura varieties introduced in the previous two sections, as well as the unitary Shimura varieties introduced in \Cref{s:unitary-sv}, are examples of Shimura varieties of \emph{PEL type}. In this section, we formulate a general modularity conjecture for special cycles on these Shimura varieties that recovers many of the previous conjectures in Types I, II, and III. The discussion in this section is inspired from Kai--Wen Lan's notes \cite{kai-wen-lan}.

\begin{definition}\label{def:pel-data}
    A PEL datum (Polarization, Endomorphism and Level structure) \[(\cO,L,*,(\cdot,\cdot),h_0)\] is given by: 
    \begin{enumerate}
    \item An order $\cO$ in a finite dimensional semisimple $\Q$-algebra $B$;
    \item A positive involution $*$ on $\cO$, i.e., an anti-automorphism of order $2$ that satisfies: 
    \[\Trace_{\cO\otimes_\Z \R/\R}(aa^*)>0,\,\forall a\in\cO\otimes_\Z \R,\,a\neq 0~.\]
    
    \item A finite rank free $\Z$-module $L$ with an action of $\cO$;
        \item A pairing \[(\cdot,\cdot):L\times L\rightarrow \Z\]
        compatible with the involution $*$: 
        \[(bx,y)=(x,b^*y)\] for all $x,y\in L$ and $b\in \cO$.
    \item $h_0$ is a morphism of $\R$-algebras:
    \[h_0: \C\rightarrow \mathrm{End}_{\cO_\R}(L_\R)\] satisfying the following properties:
    \begin{enumerate}
        \item For all $z\in\C$, $x,y\in L_\R$, we have \[(h_0(z)x,y)=(x,h_0(\overline z)\cdot y)~.\]
        \item The $\R$-bilinear pairing $(\cdot, h_0(i)(\cdot))$ is symmetric positive definite.  
    \end{enumerate}
    \end{enumerate}
\end{definition}

A PEL datum yields an abelian variety which 
 is defined as the real torus $A=(L\otimes \R)/L$ with complex structure given by $h_0$. The polarization is explicitly constructed as follows: the dual lattice $L^\vee$ to $L$ can be identified with the lattice: 
 
\[L^\vee=\{x\in L_\Q\,|\,\forall y\in L, (x,y)\in \Z\}~.\] 

Clearly $L\subset L^\vee$, which yields an isogeny of complex tori $A\rightarrow A^\vee=L_\R/L^\vee$, therefore defining a polarization on $A$.
\medskip

The action of the order $\cO$ on $L$ yields a morphism:
 \[\cO\rightarrow \mathrm{End}_{\C}(A)\]
satisfying a compatibility condition between the Rosati involution on $A$ and the involution $*$ on $\cO$ that we do not state but refer to \cite[Equation (5.1.1.9)]{kai-wen-lan}.  

\medskip

Given a PEL datum as in Definition \ref{def:pel-data}, we can define a moduli space over $\C$ parameterizing abelian varieties with PEL structure described as above, denoted $\cM$. This moduli space is closely related to certain class of Shimura varieties that we now define, see also \cite[Section 5.1.3]{kai-wen-lan}.  

Let $(\cO,L,*,(\cdot,\cdot),h_0)$ be a PEL datum. Let $\bfG$ be the algebraic group over $\Q$ such that for each $\Q$-algebra $R$, we have
\[\bfG(R)=\{g\in \GL_\cO(L\otimes _\Z R)\}|, \forall x,y\in L, (gx,gy)=\nu(g)(x,y),\, \nu(g)\in R^{\times}\}~.\]

The data $h_0:\C\rightarrow \mathrm{End}_{\cO_\R}(L\otimes_\Z \R)$ provides a morphism: 
\[h_0:\C^\times\rightarrow \bfG(\R)\]
and its $\bfG(\R)$ conjugacy class defines a manifold $\cD=\bfG(\R)\cdot h_0~.$
\medskip 

It is important to note that $(\bfG,\cD)$ does not always define a Shimura datum, but when it does, it is closely related to the moduli space $\cM$. In particular, this construction recovers all examples of Shimura varieties introduced in \cref{s:siegel-sv,s:unitary-sv,s:quaternion-SV}, as explained in \cite[Section 5.1.3]{kai-wen-lan}.

We now define special cycles in $\cM$. Let $(\cO,*,L_1,(\cdot,\cdot))$ and $(\cO,*,L_2,(\cdot,\cdot))$ be two sub-PEL datum such that $L_1\oplus L_2$ is orthogonal with respect to the pairing $(\cdot,\cdot)$ and finite index in $L$.

Let $\cM_1$ and $\cM_2$ be the PEL Shimura varieties constructed from the PEL data above. Then we obtain a finite morphism of Shimura varieties: 

\[\iota_{L_1,L_2}:\cM_1\times\cM_2\rightarrow \cM~.\] 

Let $m\geq 1$ and let $\cZ(m)$ denote the union of images of the morphisms $\iota_{L_1,L_2}$ across all PEL sub-data such that $L_1\oplus L_2$ has index $m$ in $L$.

\begin{remark}
    The cycles $\cZ(m)$ could also be defined using the modular interpretation: it is the locus of abelian varieties $A$ that are isogenous to a product $A_1\times A_2$ where $A_1$ is in $\cM_1$ and $A_2$ is in $\cM_2$.
\end{remark}

Finally, define the following generating series: 

\[\Theta_\cM(q)=\sum_{m\geq 0}[\cZ(m)]q^m\in \CH^*(\cM)\llbracket q\rrbracket~.\]

\begin{question}
    What can be said about the modularity of the generating series $\Theta_\cM$?
\end{question}

We present now some heuristics suggesting that the question above has a reasonable answer. These are based on the recent paper \cite{greer-lian}. 
\medskip

Let $\cM_1$ be the PEL moduli space of abelian varieties with PEL datum $(\cO,*,L_1,(\cdot,\cdot))$. The tensor product $W=L_1^\vee\otimes_\cO L_2=\mathrm{Mor}_\cO(L_1,L_2)$ admits a symmetric bilinear pairing $Q=\psi_1\otimes \psi_2$. Let $Y$ be the arithmetic quotient of the period domain arising from Hodge structures of weight $0$ on $W$, polarized by $Q$. We then get a period map: 

\[\cM_1\times \cM\rightarrow \Gamma\backslash\cD~.\]

The special cycles $\cZ(m)$ in $\Gamma\backslash\cD$ pull back to special cycles on $\cM_{1}\times \cM$, and their projections to $\cM$ coincide \emph{set-theoretically} with the special cycles defined above. By Kudla--Millson, the generating series of the special cycle classes $[\cZ(m)]$ are modular in the cohomology of $\Gamma\backslash \cD$, so it is natural to expect that the resulting generating series in $\cM$ is also modular. However, there are difficulties from improper intersections in general. Furthermore, since $\cM_1$ may not be compact, there is no well-defined pushforward map on cycles from $\cM_1\times \cM$ to $\cM$.
 
\begin{example}
    When $\cO=\Z$, then the Shimura datum recovers the Siegel datum introduced the \Cref{s:siegel-sv}.
\end{example}

\begin{example}
    Let $k$ be a quadratic imaginary field and let $\cO$ be the ring of integers of $k$. Let $L$ be a $\cO$-module endowed with a Hermitian pairing of signature $(a,b)$. Then the PEL Shimura datum recovers the unitary datum from \Cref{s:unitary-sv}
\end{example}

\begin{example}
    Let $\cO$ be an order in a quaternion algebra $B$ over $\Q$ ramified at $\infty$, i.e., $\B\otimes_\Z\Q\simeq \bH$ the Hamilton quaternions. Then we recover the variety introduced in \Cref{s:quaternion-SV}.
\end{example}

Finally, we formulate an analogue of Kudla's conjecture in this general PEL-type setting.
\begin{question}\label{conjecture:kudla} 
    Does there exist a smooth compactification of $\cM\hookrightarrow \overline \cM$ and corrections to the cycle closures $\overline{\cZ}(m)$ by boundary cycles, such that the resulting generating series is also a modular form valued in $\CH^*(\overline\cM)$? 
\end{question}

\section{Modularity in compactifications}\label{s:compact}

\subsection{Satake-Baily-Borel and $L^2$-cohomology}

Let $X = \Gamma\backslash \mathcal \cD$ be a connected Shimura variety defined by the Shimura datum $(\bfG,\cD)$ and let $X\hookrightarrow X^{*}$ be its Satake-Baily-Borel compactification \cite{bailyborel}. In general, $X^{*}$ is a singular complex projective variety. We recall briefly how this compactification is constructed in general.
\medskip 

Let $\bfP\subset \bfG$ be a maximal parabolic subgroup of $\bfG$. Then the Hermitian part $\bfG_P$ of the Levi factor of $\bfP$  determines a Shimura datum $(\bfG_P,\cD_P)$. It has an associated Shimura variety $X_P$. There are only finitely many $\Gamma$-equivalence classes of maximal rational parabolic subgroups $\bfP$ and the union 
\[X^{*}=X\cup\bigcup_{\bfP/\Gamma} X'_{\bfP}\]
admits the structure of a projective algebraic variety. 
\medskip 

As singular varieties, the Baily-Borel compactifications carry intersection cohomology groups, defined in \cite{goresky-macpherson}. Assume that we are given a collection of special cycles $\cZ(m)$ in $X$ of codimension $c\geq 1$, such that their generating series is a modular form valued in cohomology. The next most natural generating series to consider is the series of Zariski closures in $X^{*}$. By a well-known theorem of Looijenga \cite{Looijenga-Rapoport,Looijenga-zucker-1} and Saper--Stern \cite{saper-stern}, answering a conjecture of Zucker \cite{zucker-conjecture}, the intersection cohomology of $X^{*}$ is isomorphic to the $L^2$-cohomology of $X$. It is therefore tempting to make the following conjecture:

\begin{conjecture}
   The classes $[\cZ(m)]$ are in fact defined in the $L^2$-cohomology of $X$, and their generating series is a modular form.
\end{conjecture}

To our knowledge, the main evidence for this conjecture is \cite{bruinier-zuffetti} where the authors prove the result for orthogonal Shimura varieties in codimension up to the middle degree, and the earlier \cite{getz-goresky} in the case of Hilbert modular surfaces.

\subsection{The Boundary Stratification}

Each maximal parabolic $\bfP$ is the stabilizer of some $\Q$-isotropic subspace $W\subset V$ of dimension $r>0$. This is the {\it rank} of the associated boundary stratum. The Levi quotient of $\bfP$ has real points $G_r\times H_r$, where $G_r$ is of Hermitian type. This $G_r$ gives rise to a Shimura variety $X'_{\bf P}$, which forms a stratum in the boundary of $X^*$. When $r=1$, $H_1$ is a finite group. When $r$ takes the maximal possible value, namely the rank of $\bfG$, then the Levi quotient $H_{\max}$ has trivial Hermitian type factor, and the associated boundary stratum is a point. In that case, $H_{\max}$ acts naturally on a certain convex cone, denoted by $\mathscr C_{\max}$, which will play an important role in the toroidal compactification.

\begin{table}[!h]
\begin{center}
\begin{tabular}{|c|c|c|c|c|c|} 
 \hline
Type & $G$ & $G_1$ & $G_2\times H_2$ &  $H_{\max}$ & $\mathscr C_{\max}$  \\ \hline 
I & $\SU(a,b)$ & $\SU(a-1,b-1)$  & $\SU(a-2,b-2)\times \SL_2(\C)$ & $\SL_a(\C)$ & $\Herm^+_a(\C)$ \\ \hline
II & $\SO^*(2n)$ & $\SO^*(2n-4)$ & $\SO^*(2n-8)\times \SL_2(\H)$ & $\SL_n(\H)$ & $\Herm^+_n(\H)$ \\ \hline
III & $\Sp(2n,\R)$ & $\Sp(2n-2,\R)$ & $\Sp(2n-4,\R)\times\SL_2(\R)$ & $\SL_n(\R)$ & $\Sym^+_n(\R)$ \\ \hline
IV & $\SO(n,2)$ & $\SL(2,\R)$ & $\{1\}\times\SO(n-1,1)$ & $\SO(n-1,1)$ &  $\mathscr C^-\subset \R^{n-1,1}$ \\
\hline
\end{tabular}\bigskip
\caption{Some data associated to the boundary strata in $X^*$}
\label{table-data-boundary}
\end{center}
\end{table}

\subsection{Toroidal compactification}

A particularly nice class of smooth compactifications of Hermitian locally symmetric spaces was constructed in \cite{amrt}; these are examples of \emph{toroidal compactifications}, and they depend on a choice of combinatorial data $\Sigma$ called a fan. As smooth projective varieties, they form natural candidates for Kudla's Conjecture \ref{conjecture:kudla}.
\medskip 

The combinatorial input $\Sigma$ is a {\it rational polyhedral cone decomposition} of $\mathscr C$, for each maximal parabolic $\bfP$ that is admissible with respect to the action of $\Gamma\cap \bfP_\Q$. The resulting projective variety $X^\Sigma$ admits a birational morphism 
$$\eps:X^\Sigma \to X^*$$
which restricts to an isomorphism over $X\subset X^*$. Over each stratum $X_r'$ of rank $r$ in the boundary $X^* \smallsetminus X$, the proper morphism $\eps$ factors as:
$$\eps^{-1}(X'_r) \overset{\tau}\to \cY_r \overset{\pi}\to X'_r,$$
up to finite group actions, where $\tau$ is a fibration by toric varieties determined by the fan $\Sigma$, and $\pi$ is an abelian fibration. When $r=1$, $\tau$ will always be the identity map; we focus on this case first. In the Type IV case, $\mathcal Y_1$ is the Kuga-Sato variety over the modular curve $X_1'$. In the remaining Types I--III, $\mathcal Y_1 \to X_1'$ is the tautological PEL family of abelian varieties over $X_1'$. The dimension count is such that this gives a divisor in the toroidal boundary:

\begin{table}[!h]
\begin{center}
\begin{tabular}{|c|c|c|c|c|c|} 
 \hline
Type & $G$ & $G_1$ & $\dim(X'_1)$ & $\dim(Y_1)$ &  $\dim(\eps^{-1}(X'_1))$  \\ \hline 
I & $\SU(a,b)$ & $\SU(a-1,b-1)$  & $ab-a-b+1$ & $a+b-2$ & $ab-1$ \\ \hline
II & $\SO^*(2n)$ & $\SO^*(2n-4)$ & $(n-2)(n-3)/2$ & $2n-4$ & $n(n-1)/2 -1$ \\ \hline
III & $\Sp(2n,\R)$ & $\Sp(2n-2,\R)$ & $n(n-1)/2$ & $n-1$ & $n(n+1)/2-1$ \\ \hline
\end{tabular}\bigskip
\caption{Rank 1 boundary dimension count for PEL}
\label{table-rk1-boundary}
\end{center}
\end{table}

For rank $r>1$, there is a non-trivial toric part of the fibration, and the abelian fibration $\pi:\mathcal Y_r \to X'_r$ is a self fibered product of the tautological PEL family. More precisely, for each maximal parabolic $\bfP$ there is a rational isotropic subspace $W\subset V$, with the appropriate complex or quaternionic structure, and we take $\mathcal Y_r=\mathcal Y_1\otimes_F W$ where $F=\R$, $\C$, or $\H$ depending on the case. The fibers of $\tau:\eps^{-1}(X')\to \mathcal Y_r$ have dimension $\dim_\R(\mathscr C_r)-1$, and finite unions of toric varieties constructed from the rays of $\Sigma$.

\begin{table}[!h]
\begin{center}
\begin{tabular}{|c|c|c|c|c|c|} 
 \hline
Type & $G$ & $G_r$ & $\dim(X'_r)$ & $\dim(Y_r)$ & $\dim(\mathscr C_r)$   \\ \hline 
I & $\SU(a,b)$ & $\SU(a-r,b-r)$  & $ab-r(a+b)+r^2$ & $r(a+b-2r)$ & $r^2$  \\ \hline
II & $\SO^*(2n)$ & $\SO^*(2n-4r)$ & $(n-2r)(n-2r-1)/2$ & $r(2n-4r)$ & $r(2r-1)$  \\ \hline
III & $\Sp(2n,\R)$ & $\Sp(2n-2r,\R)$ & $(n-r) (n-r+1)/2$ & $r(n-r)$ & $r(r+1)/2$\\ \hline
\end{tabular}\bigskip
\caption{Rank $r$ boundary dimension count for PEL}
\label{table-rkr-boundary}
\end{center}
\end{table}
One can easily check that the sum of the last three columns is equal to $\dim(X)$, so the pre-image $\eps^{-1}(X'_r)$ is indeed codimension 1 inside of $X^\Sigma$. This completes our rough description of the irreducible divisor components of $X^\Sigma \smallsetminus X$.
\medskip

The limits of the special cycles $\mathcal Z(m)$ in each toroidal boundary divisor can be described as follows. In Type IV, the closure $\ov{\mathcal Z}(m)$ meets each boundary divisor of $X^\Sigma$ in a hypersurface which dominates $X'_r$ for $r=1,2$. For $r=1$, it is a union of families of abelian hypersurfaces in the Kuga-Sato variety over $X'_1$, and for $r=2$, it is a union of toric character hypersurfaces over the point $X'_2$. Types I--III are of a different nature; the closure $\ov{\cZ}(m)$ meets each boundary divisor in a codimension $c$ subvariety which dominates a special cycle $\mathcal Z'(m)\subset X'_r$. Its fibers over points in $X'_r$ are unions of (closures of) codimension $r$ semiabelian subvarieties of the fibers of $\eps$. The enumeration of semiabelian subvarieties gives rise to weighted theta series, which we expect are quasi-modular when $r=1$ and mixed mock modular when $r>1$, due to the presence of toric variety factors.
\medskip

In this next section, we make this geometry more explicit in the two examples $\bfG=\SO(n,2),\SU(n,1)$.

\subsubsection{Orthogonal case}
In this section, we summarize the main results of \cite{engel-greer-tayou}, which were obtained independently by \cite{garcia-modularity}, as well as in the earlier work \cite{bruinierzemel}.
\medskip 

We will reuse the setup from \Cref{s:orthogonal}. The Baily-Borel compactification of an orthogonal Shimura variety $X$ contains two types of boundary components ($r=1,2$), which are parametrized by $\Gamma$-equivalence classes of isotropic subspaces of $L$. The nomenclature for these types (see below) follows the conventions of Kulikov for K3 surface degenerations. The reader may prefer to use $r=1$ for Type II and $r=2$ for Type III, in order to avoid confusion with the nomenclature for Shimura variety types.

\begin{itemize}
    \item {\bf Type II components:} these are parametrized by totally isotropic primitive planes $J\subset L$ and the corresponding boundary component in $X^{*}$ is a modular curve $X'_J$.
    \item {\bf Type III components:} these are parametrized by isotropic primitive lines $I\subset L$ and the corresponding boundary component in $X^{*}$ is a point $\{*_I\}$.
\end{itemize}

As for the toroidal compactifications, they depend on the following data: for each isotropic primitive line $I$, the lattice $K=I^\bot/I$ has signature $(n-1,1)$. The fan data $\Sigma$ corresponds to a choice of an admissible cone decomposition of \[\mathscr C^-=\{x\in K_\R\,|\,Q(x)<0\}~,\]
which is invariant under the arithmetic group $\Gamma_K\subset \O(K,Q)$.
By choosing $\Sigma$ simplicial and unimodular, we can guarantee that the compactification $X^\Sigma$ is smooth. The boundary divisors are described as follows: 

\begin{itemize}
\item {\bf Type II divisors:} these are pre-images of the Type II boundary components of $X^{*}$, and they are isomoprhic to Kuga--Sato varieties $\cY\otimes_\Z J^\bot/J\rightarrow X'_J$, where $\cY\rightarrow X'_J$ is the universal elliptic curve. We denote such boundary divisors by $\cB_J$.

\item {\bf Type III divisors:} these are pre-images of the Type III boundary components of $X^{*}$, and are in bijection with $\Gamma_K$-orbits of one-dimensional rays $\sigma$ in the fan decomposition $\Sigma$ of $\mathcal C^-$. We denote the corresponding boundary divisors by $\cB_{I,\sigma}$.
\end{itemize}

Let $\theta(\tau)$ denote the usual theta series of a definite lattice, let $E_2(\tau)=1-24\underset{n\geq 1}\sum \sigma_1(n)q^n$,
and let $F^+_N(\tau)$ be the Zagier mock modular form of weight $3/2$ introduced in \cite[2.3]{engel-greer-tayou}. The following theorem settles Kudla's conjecture for cycles of codimension $1$. 

\begin{theorem}[\cite{engel-greer-tayou}]
   Assume that $n\geq 3$. The corrected generating series 
\begin{align*} 
\sum_{m\geq 0} &[\ov{\cZ}(m)]q^{m} -\frac{1}{2}\sum_{(I,\sigma)} \left(\theta_{\sigma^\perp}\cdot F^+_{N(\sigma)}\right)\otimes \cB_{I,\sigma} \\ &+ \frac{1}{12}\sum_{J}
\left(E_2\cdot \theta_{J^\bot/J}\right)\otimes \cB_J
\end{align*}
valued in $\Pic_\Q(X^\Sigma)\llbracket q \rrbracket$ is a weakly holomorphic modular form of weight $1+\frac{n}{2}$ for $\mathrm{SL}_2(\Z)$. 
\end{theorem}

Finally, we mention that Kudla's conjecture in the top codimension for orthogonal Shimura varieties has also been solved in \cite{bruinier-rosu-zemel}.

\subsubsection{Unitary case}
In this section, we summarize the main results of  \cite{greer-tayou}. We will reuse the setup from \Cref{s:unitary-sv}, with $a=n+1$ and $b=1$ and $\Gamma$ neat for simplicity. The unitary Shimura variety $X$ is a ball quotient. The Baily-Borel compactification $X^{*}$ of $X$ is obtained by adding finitely many points that correspond to $\Gamma$-equivalence classes of primitive isotropic $\cO_k$-lines $I\subset L$. Given such a line $I$, the $\cO_k$ lattice $I^\bot/I$ is endowed with a Hermitian form which positive definite. The boundary of $X$ in its canonical toroidal compactification $X^{\rm tor}$ is a disjoint union of divisors $\cB_I$ indexed by the lines $I$. Since $\Gamma$ is neat, we have $\cB_I\simeq E\otimes_{\cO_k} I^\bot/I$ where $E=\C/\cO_k$.
\medskip 

The first tool that we use is the following splitting result, which is proved using the Hard Lefschetz Theorem on each boundary divisor $\cB_I$ separately.
\begin{lemma}
    For $2g\leq n$, the following pushforward map on homology is surjective: 
    \[H_{2g}(X,\Q)\oplus\bigoplus_{\{I\}/ \Gamma} H_{2g}(\cB_I,\Q)\rightarrow H_{2g}(X^{\rm tor},\Q)~.\]
\end{lemma}

Let $2g\leq n$, and let $C\in H_{2g}(X^{\rm tor},\Q)$ be a homology class. By the lemma above, we can write $C=C_1+C_2$ where $C_1\in H_{2g}(X,\Q)$ and $C_2\in \bigoplus_{\{I\}/ \Gamma} H_{2g}(\cB_I,\Q)$.

 By the theorem of Kudla and Millson \cite{kudla-millson}, the generating series 
 \[\sum_{M\in \mathrm{Herm}_g(\cO_k)_{\geq 0}}([\cZ(M)] \cap C_1)\,q^M\]
is a Hermitian holomorphic modular form of genus $g$ and weight $n+1$ with respect to an arithmetic subgroup of $\U(g,g)(\Z)$. It is enough therefore to analyze the intersection of the generating series of special cycles with $C_2$, which takes place entirely within the boundary.
\medskip

Let $J=I^\bot/I$. The intersection of the special cycles with the boundary can be described as follows: for each $\lambda\in J$, let $Z(\lambda)\subset \cB_I$ denote the kernel of the map of abelian varieties: 
\[E\otimes_{\cO_k}J\xrightarrow{(\cdot,\lambda)} E,\]
and we assign a cycle multiplicity equal to the divisibility of $\lambda \in J$ squared.
\begin{lemma}
    The intersections of the special cycles $\cZ(M)$ with $\cB_I$ is given set-theoretically by the following sum: 
    \[\sum_{\underset{h(\underline \lambda)=M}{\underline\lambda\in J^g}} Z(\lambda_1) \cap \cdots\cap Z(\lambda_g)~.\]
\end{lemma}

According to the lemma above, the restriction of the generating series $\Theta^g_L(q)$ to $\mathcal B_I$ becomes simply the generating series: 

\[\sum_{\underline \lambda\in J^g}[Z(\underline \lambda)]q^{h(\underline \lambda)}\in\CH^g(\cB_I) \llbracket q \rrbracket~.\]

\begin{proposition}
    There exists a polynomial $f:J^g_\R\rightarrow \cH^{g,g}(E\otimes J,\C)$
such that for every $\underline \lambda\in J$, the harmonic differential form $f(\underline\lambda)$ has cohomology class equal to $[Z(\underline \lambda)]$. The polynomial $f$ satisfies the following invariance property:
\[
\forall\underline \lambda\in J_\R^g, \forall A\in \GL_g(\C),\,f(\underline \lambda\cdot A)=|\det A|^2f(\underline\lambda) ~.
\]
\end{proposition} 

The space of all polynomials $P:J_\R^g\rightarrow \C $ that satisfy this invariance condition is denoted by $\cF_{n,g}$, and it is a finite-dimensional vector space. We prove the following crucial theorem. 
\begin{theorem}\label{t:sl2-action}
    The direct sum $\bigoplus_{\ell=0}^n\cF_{n,\ell}$ admits an $\mathfrak{sl}_2$-action where the lowering operator is the Laplacian on $J_\R$, which is invariant under the action of $\U(h)$.
\end{theorem}

The conormal bundle $\cN_{\cB_I}^\vee$ of the boundary $\cB_I$ is ample and it coincides, up to a scalar, with the polarization on $E\otimes_{\cO_k} L$ constructed out of the natural polarization $E$ and the Hermitian form $h$ on $J$. Let us denote by $D\in \Pic(\cB_I)$ the corresponding divisor class. 
\medskip 

The cycle $Z(\underline\lambda)$ naturally lives in $H^{g,g}(E,\Q)\otimes \cF_{n,g}$. Using a Lefschetz-type decomposition, we are able to prove the following theorem: 

\begin{theorem}
   For each $0\leq k\leq g$, there exists a family $(W^k_i)_{i\in I_k}$ of algebraic cycles of codimension $k$  and a family of harmonic polynomials $P_i^k\in \cF_{n,k}$ uniquely determined by these cycles such that for every $\underline\lambda\in J_\R^g$:
    \[[Z(\underline \lambda)]=\sum_{\ell=0}^{g}\sum_{i\in I_{\ell}}(\Lambda^{g-\ell}P^{\ell}_i)(\underline \lambda) W^{\ell}_i\cap D^{g-\ell}~\in H^{2g}(\cB_I,\R),\]
    where $\Lambda$ is the raising operator given in \Cref{t:sl2-action}.
\end{theorem}

It follows from the theorem above that \[\cZ(\underline \lambda)-\sum_{\ell=0}^{g-1}\sum_{i\in I_{\ell}}(\Lambda^{g-\ell}P^{\ell}_i)(\underline \lambda) W^{\ell}_i\otimes D^{g-\ell}\] is a harmonic polynomial in $\underline\lambda$ and satisfies the linear invariance properties of elements in $\cF_{n,g}$. The result now follows using a theorem of Freitag and Shimura. We refer the reader to \cite[Section 5]{greer-tayou} for the end of the proof.

\section{Non-classical period domains}\label{s:vhs}

In this section, we speculate on modularity theorems for more general period domains. 
\subsection{Variations of Hodge structure}
Let $\bV=\{\bV_\Z,\cF^\bullet\cV,Q,\nabla\}$ be a pure polarized variation of Hodge structure of weight $2k$ over a smooth, connected, quasi-projective algebraic variety $\cS$. It consists of the data of a local system $\bV_\Z$ of finite rank free $\Z$-modules, together with a quadratic form $Q$ on $\bV_\Z$, and a holomorphic filtration $\cF^{2k}\subset\ldots\subset \cF^0=\cV$ on the flat bundle $(\cV=\bV_\Z\otimes \cO_S,\nabla)$ associated to $\bV_\Z$ by the Riemann--Hilbert correspondence. 
\medskip 

Examples of such data are provided by variations of \emph{geometric origin}: if $f:\cX\rightarrow \cS$ is a smooth projective morphism, then $R^{2k}f_*\Z/\textrm{torsion}$ is a local system on $\cS$. The flat vector bundle $R^{2k}f_*\Z\otimes_\C \cO_\cS\simeq R^{2k}f_*(\Omega^{\bullet}_{\cX/\cS})$ admits a Hodge filtration induced by the stupid filtration on the hypercohomology complex $R^{2k}f_*\Omega^{\bullet}_{\cX/\cS}$. Finally, the quadratic form is obtained using the Poincar\'e pairing, using a $f$-relatively ample line bundle. 

\medskip 

We restrict our attention in what follows to the case of weight $2$ pure variations. Let $\pi:\widetilde{\cS} \to \cS$ be the universal cover of $\cS$, and let $L=H^0(\widetilde S,\pi^*\bV_\Z)$. Then $(L,Q)$ is a free $\Z$-module endowed with a quadratic form of signature $(n,2p)$. The period domain associated to this variation is equal to:
\begin{align*}
    \cD&=\{W\subset L_\C\,\big| \, \dim(W)=p,\, Q|_W=0,\, Q(v,\overline v)>0,\, \forall v\in W\smallsetminus \{0\}\}\\
    &\simeq \SO(n,2p)/\SO(p)\times \U(p)~.
\end{align*}    
For $p>1$, $\cD$ is {\it not} the symmetric space of $\SO(n,2p)$, the latter being equal to \[\cD_{\rm sym}=\SO(n,2p)/S(\O(n)\times \O(2p))~.\]
Let $\Gamma=O(L)$, or any congruence subgroup containing the image of the monodromy group of the local system $\mathbb V_\Z$.
The complex manifold $\Gamma\backslash \cD$ admits a family of special cycles which were defined by Kudla and Millson: 
let $\lambda\in L^\vee$ such that $Q(\lambda)>0$. Then the subvariety 
\[\cD(\lambda)=\{W\in\cD\,|\, \lambda \textrm{ is orthogonal to } W\}\]
is a sub-domain associated to the group $\SO(n-1,2p)$.

Let $m\in \Q_{>0}$, and set 
$\cZ(m)=\Gamma\backslash \left(\underset{Q(\lambda)=m}\bigcup \cD(\lambda)\right)$. It follows from the Kudla--Millson theorem \cite{kudla-millson}, pulled back from the locally symmetric space $\Gamma\backslash\cD_{\rm sym}$, that the generating series: 
\[\Phi(q)=\sum_{m\geq 0}[\cZ(m)]q^m\in H^{2p}(\Gamma\backslash \cD,\Q) \llbracket q \rrbracket\] is the $q$-expansion of a holomorphic modular form. By further pullback along the classifying morphism $\mathcal S \to \Gamma\backslash \cD$, the generating series of Noether--Lefschetz loci:

\[
\Phi^{\NL}(q)=\sum_{m\geq 0}[\NL(m)]q^m\in H^{2p}(\cS,\Q)\llbracket q \rrbracket
\]
is also a holomorphic modular form of weight $(n+2p)/2$. Since $\cS$ is an algebraic variety, it is natural to ask the following question.

\begin{question}
    Is the Noether--Lefschetz generating series a modular form valued in $\CH^p(\cS)$?
\end{question}


\subsection{Compactifications}

The complex varieties $\Gamma\backslash \cD$ are non-algebraic as soon as $p>1$; see \cite{griffiths-robles}. It is therefore meaningless to formulate Kudla's Conjecture \ref{conjecture:kudla} for these period domains quotients. However, the algebraic variety $\cS$ admits projective compactifications, and it is still interesting to study the modularity behavior of Noether--Lefschetz cycle closures there. 
\medskip 

Consider a simple normal crossings projective compactification: 
\[\cS\hookrightarrow \overline \cS,\]
and assume that all the local monodromy operators are unipotent.  

\begin{question}
    There exist suitable corrections $\widetilde{\NL}(m)$ of the Zariski closures $\overline{\NL}(m)$ by cycles supported in the boundary $\overline{\cS} \smallsetminus \cS$, such that their resulting generating series is a holomorphic modular form with values in $\CH^{p}(\overline\cS)_\Q$.
\end{question}

This level of generality is amenable to the methods in \cite{garcia-modularity}, where the Kudla-Millson form pulled back to $\mathcal S$ is shown to be integrable when $p=\dim(\mathcal S)=1$.


\section{Curious coincidences}\label{s:coincidence}

There is a short list of exceptional isogenies between certain pairs of small-dimensional groups from disparate families in Table \ref{table_sv}. We record them here for reference:

\begin{minipage}{0.48\textwidth}
\begin{align}
\mathrm{SU}(1,1) &\to \mathrm{SO}(1,2)\\
\mathrm{Sp}(4,\mathbb{R}) &\to \mathrm{SO}(3,2)\\
\mathrm{SU}(2,2) &\to \mathrm{SO}(4,2)
\end{align}
\end{minipage}
\hfill 
\begin{minipage}{0.48\textwidth}
\begin{align}
\mathrm{SU}(3,1) &\to \mathrm{SO}^*(6) \\
\mathrm{SO}^*(8) &\to \mathrm{SO}(6,2)
\end{align}
\end{minipage}
\medskip

The groups in each pair have the same Hermitian symmetric domain $\cD$, so in particular they can give rise to the same Shimura variety $X$. At the same time, the codimensions and weights of the generating functions for special cycles can look quite different, depending on which group from the isogenous pair is used.

\begin{question}
    For each isogeny $G_1 \to G_2$ from the list above, we can arrange two associated generating series $\Phi_i(q)$, $i=1,2$, valued in $\CH^d(X)$ for the {\it same} Shimura variety $X$:
    $$\Phi_i(q) \in \CH^{d}(X)\otimes_\Q \Mod_{k_i}(\Sp_{2g_i}(\Z)) ~.$$
Are $\Phi_1$ and $\Phi_2$ related by a natural lifting between spaces of modular forms? 
\end{question}

For the isogenies (1) and (2), the two different generating series of special divisors are likely related by a version of the Shimura lift \cite{shimuralift}: 
$$\xi_{\rm{Sh}}: \Mod_{k+1/2}(\SL_2(\Z)) \to \Mod_{2k}(\SL_2(\Z)).$$
This perspective for isogeny (1) has been explained in detail by Sankaran \cite{sanka1},\cite{sanka2} for arithmetic Chow groups of modular curves. 

For the isogeny (3), the expected relation is via the Saito-Kurokawa lift \cite{zagsk}:
$$\xi_{\mathrm{SK}}: \Mod_{2k}(\SL_2(\Z)) \to \Mod_{k+1}(\Sp_4(\Z)).$$

\begin{conjecture}
Let $X$ be a Shimura variety for $\SU(2,2)\sim \SO(4,2)$ that is simultaneously unitary and orthogonal type. The generating series $\Phi_2(q)\in \CH^2(X)\otimes \Mod_3(\Sp_4(\Z))$ of codimension 2 orthogonal special cycles is the Saito-Kurokawa lift of the generating series $\Phi_1(q)\in \CH^2(X) \otimes \Mod_4(\SL_2(\Z))$ of codimension 2 unitary special cycles.
\end{conjecture}

The isogeny (5) was recently given a modular interpretation in \cite{poon}, using the Kuga-Satake construction. We expect an Ikeda lift \cite{ikeda} to relate pairs of generating series of codimension 4 special cycles on the associated Shimura varieties. The Hermitian version of the Ikeda lift \cite{ikeda2} should relate pairs of generating series for codimension 2 special cycles in the case of isogeny (4).

Going one step further, it is natural to ask whether the same lifting relations hold between pairs of {\it Kudla-corrected} generating series with values in $\CH^d(X^\Sigma)$, after a suitable choice of smooth toroidal compactification.

\bibliographystyle{alpha}
\bibliography{bibliographie}

\end{document}